\documentclass[12pt,oneside]{article}

\usepackage{amsthm,amsfonts,mathrsfs,  amsbsy, amssymb,amsmath,graphicx,euscript}

\usepackage{geometry}
\theoremstyle{definition}
\newtheorem{definition}{Definition}

\newtheorem{remark}{Remark}

\newtheorem{ass}{Assumption}
\usepackage[english, russian]{babel}

\theoremstyle{plain}
\newtheorem{theorem}{Theorem}

\begin{document}

\begin{center} {\bf Chernoff approximations as a method for finding the resolvent of a linear operator and solving a linear ODE with variable coefficients\\

Ivan\,D.\,Remizov$^{1,2}$}

$^{1}$ Laboratory of Topological Methods in Dynamics, NRU HSE

25/12 Bol. Pecherskaya Ulitsa, Room 224, Nizhny Novgorod, 603155, Russia

$^{2}$ Institute for Information Transmission Problems of the Russian Academy of Sciences

Bolshoy Karetny per. 19, build.1, Moscow 127051 Russia

ivremizov@yandex.ru
\end{center}

\textbf{Abstract.} The Chernoff approximation method is a powerful and flexible tool of functional analysis, which allows in many cases to express exp(tL) in terms of variable coefficients of a linear differential operator L. In this paper, we prove a theorem that allows us to apply this method to find the resolvent of L. Our theorem states that the Laplace transforms of Chernoff approximations of a $C_0$-semigroup converge to the resolvent of the generator of this semigroup. We demonstrate the proposed method on a second-order differential operator with variable coefficients. As a consequence, we obtain a new representation of the solution of a nonhomogeneous linear ordinary differential equation of the second order in terms of functions that are coefficients of this equation, playing the role of parameters of the problem. For the Chernoff function, based on the shift operator, we give an estimate for the rate of convergence of approximations to the solution.\\



\begin{center} {\bf Черновские аппроксимации как метод нахождения резольвенты линейного оператора и решения линейного ОДУ с переменными коэффициентами\\
		
И\,Д.\,Ремизов$^{1,2}$}

	$^{1}$Лаборатория топологических методов в динамике НИУ ВШЭ
	
	603155, Россия, Нижний Новгород, Большая Печерская, 25/12
	
	$^{2}$Институт проблем передачи информации имени А. А. Харкевича РАН
	
127051, г. Москва, Большой Каретный переулок, д.19 стр. 1.
	
	ivremizov@yandex.ru
\end{center}

\textbf{Аннотация.} Метод черновских аппроксимаций является мощным и гибким инструментом функционального анализа, позволяющим во многих случаях выразить exp(tL) через переменные коэффициенты линейного дифференциального оператора L. В данной работе доказывается теорема, позволяющая применять этот метод для нахождения резольвенты оператора L. Наша теорема утверждает, что преобразования Лапласа аппроксимаций Чернова $C_0$-полугруппы сходятся к резольвенте генератора этой полугруппы. Мы демонстрируем предложенный метод на дифференциальном операторе второго порядка с переменными коэффициентами. В качестве следствия мы получаем новое представление решения неоднородного линейного обыкновенного дифференциального уравнения второго порядка в терминах функций, являющихся коэффициентами этого уравнения, играющих роль параметров задачи. Для функции Чернова на основе оператора сдвига мы даем оценку скорости сходимости приближений к решению.\\

\textbf{Ключевые слова:} полугруппы операторов, резольвента оператора, линейное ОДУ с переменными коэффициентами, представление решения, черновские аппроксимации\\

\textbf{УДК} 517.986.7+517.988.8+517.926.4

\section{Introduction}

The straightforward way to find the resolvent $(\lambda I-L)^{-1}$ of the linear operator $(L, D(L))$ in the Banach space $\mathcal{F}$ is to study (for given $\lambda\in\mathbb{C}$) if the equation $\lambda f-Lf=g$ has for each $g\in\mathcal{F}$ a unique solution $f\in D(L)$. If we know that $(L,D(L))$ is closed then this is enough, but if not then we must also check that $f$ depends on $g$ continuously. In case of success we have $f=(\lambda I-L)^{-1}g$. For example, if $\mathcal{F}$ is one of the subsets of space of all functions $f\colon\mathbb{R}\to\mathbb{C}$, and $Lf=af''+bf'+cf$ is a differential operator defined by constant coefficients $a,b,c\in\mathbb{C}$, then this task is not difficult. Indeed $\lambda f-Lf=g$ is a second order differential equation $\lambda f(x)-af''(x)-bf'(x)-cf(x)=g(x)$, $x\in\mathbb{R}$, with only one variable coefficient $g$, and there are standard formulas in ODE books to solve this equation using the variation of parameters method. But if $a,b,c$ are not constants, but functions that depend on $x$ then the situation becomes much worse because there are no standard formulas for solution of equation $\lambda f(x)-a(x)f''(x)-b(x)f'(x)-c(x)f(x)=g(x)$, $x\in\mathbb{R}$.

However, if $(L,D(L))$ holds one additional property, then we can use a non-straightforward method of finding the resolvent, this method is the main result of the paper. The property is that $(L,D(L))$ is the generator of a $C_0$-semigroup, which informally means that exponent $e^{tL}$ exists as a linear bounded operator, and depends on $t\in[0,+\infty)$ continuously in some sense. If $(L,D(L))$ is a linear bounded operator with $D(L)=\mathcal{F}$ then this condition holds automatically, but if $(L,D(L))$ is not bounded then there is a beautiful theory around it, see e.g. \cite{Ar2002,ABHN2001,EN,Gold, BS}. In most interesting cases (e.g. if $L$ is a differential operator with variable coefficients) it is difficult to calculate $e^{tL}$ directly because power series $e^{tL}=\sum_{n=0}^\infty(tL)^n/n!$ is useless in case of unbounded operator $L$. Yet it is possible to find $e^{tL}$ approximately using the Chernoff theorem \cite{Chernoff}, this theorem will be discussed below. After that, having $e^{tL}$ we calculate the resolvent via well known \cite{EN} formula $(\lambda I -L)^{-1} f=\int_0^\infty e^{-\lambda t}e^{tL}fdt$. The method may look complicated because of the many steps it has (find Chernoff function for $L$, find $e^{tL}$ via Chernoff approximations, find $(\lambda I -L)^{-1}$ via $e^{tL}$), however all these steps are supported by known methods, and so our proposed method may still be much simpler than considering the equation $\lambda f-Lf=g$. 

Brief history and overview of the results obtained up to 2017 in constructing Chernoff approximations of $e^{tL}$ for several classes of  operators $L$ can be found in \cite{Butko-2019}. Several papers on the topic showing the diversity of cases studied are \cite{R-2017,R-PotAn2020, R5, R2, RemAMC2018, RemJMP}, see also \cite{OSS2019, ST2020, KOS2023}. Speed of convergence of Chernoff approximations were studied in \cite{Drag2023SVMO, GR-2021-arxiv, GR-2021, GSK2019, VVGKR, Zag2020, Zag2022}, see also references therein. See also paper \cite{K2023} which numerically investigates the efficiency of the Monte Carlo method based on the application of the Chernoff theorem, and paper \cite{KOS2022} that mathematically substantiate such an approach to the Chernoff approximations.

Summing up, we can say that the method of Chernoff approximation is an extremely effective tool for expressing $e^{tL}$ in terms of variable coefficients of operator $L$. The present paper shows that this method can be also be used for expressing $(\lambda I-L)^{-1}$ in terms of variable coefficients of operator $L$, and for finding the solution of the corresponding  differential equation $\lambda f-Lf=g$. As an example we consider second order linear ODE with variable coefficients. 

Limits of multiple integral as multiplicity tends to infinity (such expressions are called Feynman formulas \cite{R-2017}) are one of the ways (which actually was used originally by Richard Feynman \cite{F1,F2}) to define Feynman path integral \cite{Maz,JL}. So in theorem \ref{ODE-th} a solution of an ODE is in the first time in history of science represented via Feynman formula, which can also be interpreted as Feynman integral if one wishes to do. Such theoretical step is novel and may lead to some unexpected applications in the future. 

\section{Preliminaries: operator semigroups and their Chernoff approximations}

Let us recall some relevant  definitions and facts of $C_0$-semigroup theory following \cite{EN}.

\begin{definition}\label{semigrdef} Let $\mathcal{F}$ be a Banach space over the field $\mathbb{R}$ or $\mathbb{C}$. Let $\mathscr{L}(\mathcal{F})$ be the set of all bounded linear operators in $\mathcal{F}$. Suppose we have a mapping $V\colon [0,+\infty)\to \mathscr{L}(\mathcal{F}),$ i.e.\ $V(t)$ is a bounded linear operator $V(t)\colon \mathcal{F}\to \mathcal{F}$ for each $t\geq 0.$ The mapping $V$, or equivalently the family $(V(t))_{t\geq 0}$, is called \textit{a strongly continuous one-parameter semigroup of linear bounded operators} (or just \textit{a $C_0$-semigroup}) iff it satisfies the following three conditions: 
	
	1) $V(0)$ is the identity operator $I$, i.e. $V(0)\varphi=\varphi$ for each $\varphi\in \mathcal{F}$; 
	
	2) $V$ maps the addition of numbers in $[0,+\infty)$ into the composition of operators in $\mathscr{L}(\mathcal{F})$, i.e. for all $t\geq 0$ and all $s\geq 0$ we have $V(t+s)=V(t)\circ V(s),$ where for each $\varphi\in\mathcal{F}$ the notation $(A\circ B)(\varphi)=A(B(\varphi))=AB\varphi$ is used;
	
	3) $V$ is continuous with respect to the strong operator topology in $\mathscr{L}(\mathcal{F})$, i.e. for all $\varphi\in \mathcal{F}$ vector $V(t)\varphi$ depends on $t$ continuously, i.e. the function $[0,+\infty)\ni t\longmapsto V(t)\varphi\in \mathcal{F}$ is continuous.
\end{definition}

\begin{remark}
	The definition of a \textit{$C_0$-group} $(V(t))_{t\in\mathbb{R}}$ is obtained by substituting $[0,+\infty)$ with $\mathbb{R}$ in the definition above.
\end{remark}

\begin{definition}\label{defgen}
	Let $(V(t))_{t\geq 0}$ be a $C_0$-semigroup in Banach space $\mathcal{F}$. Its \textit{infinitesimal generator} (or just \textit{generator}) is defined as the operator $L\colon D(L)\to\mathcal{F}$ with the domain 
	$$
	D(L)=\left\{\varphi\in \mathcal{F}: \textrm{there exists a limit }\lim_{t\to +0}\frac{V(t)\varphi-\varphi}{t}\right\} \subset \mathcal{F},
	$$ 
	and
	$$L\varphi=\lim_{t\to +0}\frac{V(t)\varphi-\varphi}{t}.$$ 
	
\end{definition}

\begin{remark}\label{generdef}
It is known \cite{EN} that for each $C_0$-semigroup $(V(t))_{t\geq 0}$ in Banach space $\mathcal{F}$, the set $D(L)$ is a dense linear subspace of $\mathcal{F}$. 
Moreover, $(L,D(L))$ is a closed linear operator that uniquely defines the semigroup $(V(t))_{t\geq 0}$. Also for each $C_0$-semigroup $(V(t))_{t\geq 0}$ in Banach space there exist constants $M\geq 1$ and $\omega\in\mathbb{R}$ such that $\|V(t)\|\leq Me^{\omega t}$ for all $t\geq 0$.
\end{remark}

\begin{remark}	
Very often the notation $V(t)=e^{tL}$ is used. This is a good notation for several reasons. First, properties $e^{0\cdot L}=I$ and $e^{(t+s)L}=e^{tL}e^{sL}$ are consistent with the case when $L$ is a number or a matrix. Second, if $L$ is a bounded linear operator then the exponent can be defined via the standard power series $e^{tL}=\sum_{n=0}^\infty(tL)^n/n!$ that converges in the operator norm. Finally, in general case we have $e^{tL}f=(I+tL)f+o(t)$ for all $f\in D(L)$ which again retains properties of the exponent in number and matrix case.
	
\end{remark}

Now we are ready to state the Chernoff theorem. From several options (see \cite{EN, BS, Chernoff, JL}), we choose the one given in~\cite{BS} (in equivalent formulation):
\begin{theorem}[\textsc{P.\,R.~Chernoff (1968)}, cf.~\cite{EN, BS, Chernoff, JL}]
	\label{ChernoffTheor}
	Suppose that the following three conditions are met:
	\begin{enumerate}
		\item 
		$C_0$-semigroup $(e^{tL})_{t\ge 0}$ with generator $(L,D(L))$ in Banach space $\mathcal{F}$ is given, 
		such that for some $w\geq 0$ the inequality $\|e^{tL}\| \le e^{wt}$ holds 
		for all $t\ge0$. 
		
		\item 
		There exists a strongly continuous mapping $S\colon[0,+\infty)\to\mathscr{L}(\mathcal{F})$ such that $S(0)=I$ and the inequality $\|S(t)\| \le e^{wt}$ holds for all $t\ge0$.
		
		\item 
		There exists a dense linear subspace $D\subset\mathcal{F}$ such that for all $f\in D$ there exists a limit
		$S'(0)f := \lim_{t\to +0} (S(t)f-f)/t$.
		Moreover, $S'(0)$ on $D$ has a closure that coincides with the generator $(L,D(L))$.
	\end{enumerate}
	Then the following statement holds:
	\begin{enumerate}
		\item[(C)] for every $f\in \mathcal{F}$, as $n\to\infty$ we have $S(t/n)^n f \to e^{tL}f$ locally uniformly with respect to $t\ge0$, i.e. for each $T>0$ and each $f\in \mathcal{F}$ we have
		$$
		\lim_{n\to\infty}\sup_{t\in[0,T]}\|S(t/n)^n f - e^{tL}f\| = 0.
		$$
	\end{enumerate}
Above $S(t/n)^n=\underbrace{S(t/n)\circ\dots\circ S(t/n)}_n$ is the composition of $n$ copies of linear bounded operator $S(t/n)$ defined everywhere on $\mathcal{F}$.

\end{theorem}

\begin{definition} \label{defChFun} 
	Let $C_0$-semigroup $(e^{tL})_{t\ge 0}$ with generator $L$ in Banach space $\mathcal{F}$ be given. 
	The mapping $S\colon [0,+\infty)\to \mathscr{L}(\mathcal{F})$ is called a \textit{Chernoff function for operator $L$} iff it satisfies the
	condition (C) of Chernoff theorem~\ref{ChernoffTheor}.
	In this case expressions $S(t/n)^n$ are called \emph{Chernoff approximations to the semigroup $e^{tL}$}.
\end{definition}

There is also a variant of the Chernoff theorem that allows us to prove the existence of the semigroup.

\begin{theorem}\label{FormulaChernova2}(Chernoff-type theorem, corollary 5.3 from theorem 5.2 in \cite{EN})
Let $\mathcal{F}$ be a Banach space, and $\mathscr{L}(\mathcal{F})$ be the space of all linear bounded operators on $\mathcal{F}$. Suppose there is a function
$ S\colon [0,+\infty)\to \mathscr{L}(\mathcal{F}),$
meeting the condition $S(0)=I$, where $I$ is the identity operator.
Suppose there are numbers $M\geq 1$ and $\omega\in\mathbb{R}$ such that $\|S(t)^k\|\leq Me^{k\omega t}$ for every $t\geq 0$ and every $k\in \mathbb{N}$. Suppose the limit
$\lim_{t\to +0}\frac{S(t)\varphi - \varphi}{t}=:L\varphi$
exists for every $\varphi\in \mathcal{D}\subset \mathcal{F},$ where $\mathcal{D}$ is a dense subspace of $\mathcal{F}$. Suppose there is 
a number $\lambda_0>\omega$ such that $(\lambda_0I-L)(\mathcal{D})$ is a dense subspace of $\mathcal{F}$.
	
Then the closure $\overline{L}$ of the operator $L$ is a generator of a strongly continuous semigroup 
of operators $(e^{\overline{L}})_{t\geq 0}$ given by the formula
$e^{\overline{L}}\varphi=\lim_{n\to\infty}S(t/n)^n\varphi$
where the limit exists for every  $\varphi\in \mathcal{F}$ and is uniform  with respect to $t\in [0,T]$ for every $T>0$. Moreover $(e^{\overline{L}})_{t\geq 0}$ satisfies the estimate $\|e^{\overline{L}}\|\leq M e^{\omega t}$ for every $t\geq 0$.
\end{theorem}

\begin{remark}
There are several theorems that help to find out if some linear operator (or the closure of it) generates a $C_0$-semigroup. Most general are Hille-Yosida theorem and Feller-Miyadera-Phillips theorem. Unfortunately both heavily use properties of the resolvent which is a kind of circulus vitiosus because we start all this activity to find the resolvent but need its properties to do it. 

However there are several results that only use properties of the operator to prove that it generates a $C_0$-semigroup:
\begin{itemize}
\item Stone's generation theorem: if $A$ is a self-adjoint operator in Hilbert space, then $iA$ generates a $C_0$-semigroup; even more, this is a $C_0$-group of unitary operators.

\item Lumer-Phillips theorem: linear, closed, densely defined in Banach space, dissipative operator $A$ with a property that $A - \lambda_0I$ is surjective for some $\lambda_0> 0$, generates a $C_0$-semigroup of contraction operators.

\item A. Yu. Neklyudov's inversion of Chernoff's theorem \cite{Nekl2008}.
\end{itemize}

 \end{remark}

\begin{remark}Chernoff's theorem is a  deep result of functional analysis and is designed for dealing with infinite-dimensional spaces $\mathcal{F}$. However, it can be illustrated in one-dimensional setting  in two ways, which helps to build intuition. First, one-dimensional version of Chernoff's theorem, when $\mathcal{F}=\mathscr{L}(\mathcal{F})=\mathbb{R}$, says that if $s\colon[0,+\infty)\to\mathbb{R}$, $l\in\mathbb{R}$ and $s(t)=1+tl+o(t)$ as $t\to0$, then $\lim_{n\to\infty}(s(t/n))^n=e^{tl}$, which is a simple fact of calculus. Second, one can see that Chernoff's theorem is an infinite-dimensional analogue of the forward Euler method for solving ordinary differential equations numerically.
	
\end{remark}

\section{Main result}

\begin{theorem}\label{mainth} Let $\mathcal{F}$ be real or complex Banach space, and let $\mathscr{L}(\mathcal{F})$ be the set of all linear bounded operators in $\mathcal{F}$. Suppose that linear operator $L\colon \mathcal{F}\supset D(L)\to\mathcal{F}$ generates $C_0$-semigroup $(e^{tL})_{t\geq 0}$ satisfying for some constants $M\geq 1$ and $\omega\geq 0$ inequality $\|e^{tL}\|\leq Me^{\omega t}$ for all $t\geq 0$. Suppose that function $S\colon [0,+\infty)\to \mathscr{L}(\mathcal{F})$ is given and $\|S(t)^k\|\leq Me^{\omega tk}$ for all $t\geq 0$ and all $k=1,2,3,\dots$ Let us denote the resolvent of $(L,D(L))$ by the symbol $R_\lambda=(\lambda I-L)^{-1}$ for all $\lambda\in\rho(L)$. Suppose that the number $\lambda\in\mathbb{C}$ is given and $Re\lambda>\omega$. Then $\lambda\in\rho(L)$ and:

1. If for all $T>0$ we have  $\lim_{n\to\infty}\sup_{t\in[0,T]}\left\|e^{tL}f - (S(t/n))^nf \right\|=0$ for all $f\in\mathcal{F}$, then for all $f\in\mathcal{F}$  we have
\begin{equation}\label{first}
\lim_{n\to\infty}\left\|R_\lambda f - \int_0^\infty e^{-\lambda t}(S(t/n))^nfdt \right\|=0.
\end{equation} 

2. If for all $T>0$ we have  $\lim_{n\to\infty}\sup_{t\in[0,T]}\left\|e^{tL} - (S(t/n))^n \right\|=0$,  then we have
\begin{equation}\label{second}
\lim_{n\to\infty}\left\|R_\lambda - \int_0^\infty e^{-\lambda t}(S(t/n))^ndt \right\|=0.
\end{equation}
\end{theorem}

\begin{proof}
Integrals in (\ref{first}) and (\ref{second}) can be understood as improper Riemann integrals because for each $n$ integrands are continuous functions of $t$. Moreover $\|e^{-\lambda t}(S(t/n))^n\|\leq e^{-t Re\lambda}Me^{(t/n)\omega n}=Me^{t(\omega -Re\lambda)}$ with $\omega -Re\lambda<0$, so both integrals converge. Theorem II.1.10 from \cite{EN} says that $\lambda\in\rho(L)$ and $R_\lambda f=\int_0^\infty e^{-\lambda t}e^{tL}fdt$ for each $f\in\mathcal{F}$ with $\|R_\lambda\|\leq M/(Re\lambda -\omega)$.

Let us prove item 1. Suppose $\varepsilon>0$ and $f\in\mathcal{F}$ are given. Let us prove that there exists $n_0\in\mathbb{N}$ such that for all $n>n_0$ the estimate $
\left\|R_\lambda f - \int_0^\infty e^{-\lambda t}(S(t/n))^nfdt\right\|<\varepsilon
$ holds.

Integral in (\ref{first}) and $\int_0^\infty e^{-\lambda t}e^{tL}fdt$ both converge hence the integral on the right hand side of the equality
$$
R_\lambda f - \int_0^\infty e^{-\lambda t}(S(t/n))^nfdt=\int_0^\infty e^{-\lambda t}(e^{tL}f - (S(t/n))^nf)dt
$$
converges. Moreover, it convergences uniformly in $n\in\mathbb{N}$ due to the estimate 
$$
\|e^{-\lambda t}(e^{tL} - (S(t/n))^n))\|\leq e^{-Re \lambda t}(\|e^{tL}\| + \|(S(t/n))^n\|)\leq e^{-Re \lambda t}(Me^{\omega t}+Me^{(t/n)\omega n)})\leq 2Me^{t(\omega -Re\lambda)} 
$$
with $\omega -Re\lambda<0$. 

Let us use the so-called $\varepsilon/2$-method to prove that $\int_0^\infty<\varepsilon$, using representation $\int_0^\infty=\int_0^T+\int_T^\infty$. First, using $\varepsilon$, we find such $T>0$ that $\int_T^\infty\leq\varepsilon/2$ for all $n$, and then for this $T$ we find such $n_0$ that for all $n>n_0$ we have $\int_0^T<\varepsilon/2$. This will give us $\int_0^\infty<\varepsilon$. Indeed, we have
$$
\left\|\int_T^\infty e^{-\lambda t}(e^{tL}f - (S(t/n))^nf)dt\right\|\leq \int_T^\infty 2Me^{t(\omega -Re\lambda)}dt=2M\frac{1}{Re\lambda-\omega}e^{T(\omega -Re\lambda)}\leq\varepsilon/2
$$ 
for all $n\in\mathbb{N}$ and all $T$ satisfying inequality 
$
T\geq\max\left(0,\frac{1}{Re\lambda-\omega}\ln\frac{4M}{(Re\lambda-\omega)\varepsilon}\right)$. Suppose that such number $T$ is selected.

Next, thanks to the conditions of the theorem we have $\lim_{n\to\infty}\sup_{t\in[0,T]}\left\|e^{tL}f - (S(t/n))^nf \right\|=0$, so there exists $n_0$ such that for all $n>n_0$, we have
\begin{multline*}
\left\|\int_0^T e^{-\lambda t}(e^{tL}f - (S(t/n))^nf)dt\right\|\leq \int_0^T e^{-tRe\lambda}\|e^{tL}f - (S(t/n))^nf\|dt\leq\\
\leq T\max_{t\in[0,T]}e^{-tRe\lambda}\sup_{t\in[0,T]}\left\|e^{tL}f - (S(t/n))^nf\right\|<\varepsilon/2.
\end{multline*}
So we proved that for arbitrary $\varepsilon>0$ there exists $n_0$ such that for all $n>n_0$ we have
$$
\left\|R_\lambda f - \int_0^\infty e^{-\lambda t}(S(t/n))^nfdt\right\|<\varepsilon.
$$

Item 1 is proved. Proof of item 2 is obtained by deleting $f$ from the proof of item 1.
\end{proof}

\section{Corollary: Feynman formula for the resolvent}

\begin{ass}\label{Hbar}

Consider functions $a,b,c\colon \mathbb{R}\to \mathbb{R}$ on $\mathbb{R}$. Assume that $a(x)>0$ for all $x\in \mathbb{R}$. Assume that there exists such $\beta\in(0,1]$ that function $c$ is bounded and  H\"older continuous with H\"older exponent $\beta$, and functions $a$, $x\mapsto 1/a(x)$, $b$ are bounded and H\"older continuous with H\"older exponent $\beta$ with derivatives of order one and two.

Consider Banach space $C_0(\mathbb{R},\mathbb{R})$ of all continuous functions $f\colon \mathbb{R}\to \mathbb{R}$ vanishing at infinity (i.e. $\lim_{\|x\|\to\infty}f(x)=0$), with the uniform norm $\|f\|=\sup_{x\in \mathbb{R}}|f(x)|$. Consider space $C^2_c(\mathbb{R},\mathbb{R})$ of all compactly-supported functions that are continuous with derivatives of order one and two; note that this space is a dense linear subspace in $C_0(\mathbb{R},\mathbb{R})$.

For all $f\in C^2_c(\mathbb{R},\mathbb{R})$ define linear operator $H$ via the formula
\begin{equation}\label{Hlab}
(Hf)(x)=a(x)f''(x)+b(x)f'(x)+b(x)f(x).
\end{equation}

Assume that the closure $(\overline{H},D(\overline{H}))$ of operator $H\colon D(H)= C^2_c(\mathbb{R},\mathbb{R})\to C_0(\mathbb{R},\mathbb{R})$ exists and generates a $C_0$-semigroup $(e^{t\overline{H}})_{t\geq 0}$ in $C_0(\mathbb{R},\mathbb{R})$; this assumption is fulfilled, e.g. if $a(x)\geq a_0$ for some constant $a_0>0$, and  $b,c$ only satisfy already mentioned conditions. Now operator $(\overline{H},D(\overline{H}))$ is well defined.
\end{ass}

\begin{remark}\label{remhbar}
Following \cite{BGS2010}, let us consider operator-valued functions  $S_1,S_2,S_3,S$ that are defined on $[0,\infty)$ and take values in $\mathscr{L}(C_0(\mathbb{R},\mathbb{R}))$. For all  $x,y\in\mathbb{R}, t>0$ define
$$(S_1(t)f)(x)=\frac{1}{\sqrt{4\pi ta(x)}}\int_{\mathbb{R}}\exp\left(\frac{-(x-y)^2}{4ta(x)}\right)f(y)dy,$$
$$(S_2(t)f)(x)= \frac{1}{\sqrt{4\pi ta(x)}}\int_{\mathbb{R}}\exp\left(\frac{-(x-y)^2}{4ta(x)}-\frac{b(x)(x-y)}{2a(x)}\right)f(y)dy,
$$
\begin{equation}\label{Slab}
(S_3(t)f)(x)=\exp\left(t\left(c(x)-\frac{b(x)^2}{4a(x)}\right)\right)f(x),\quad S(t)=S_3(t)S_2(t).
\end{equation}

It is proved in \cite{BGS2010} that for all $f\in C^2_c(\mathbb{R},\mathbb{R})$ and $t\to 0$ we have
$$(S_1(t)f)(x)=f(x)+ta(x)f''(x)+o(t),$$
$$ (S(t)f)(x)=f(x) +t[a(x)f''(x)+b(x)f'(x)+c(x)f(x)]+o(t),
$$
and inequality $\|S(t)\|\leq e^{w\cdot t}$ holds for all $t\geq 0$ with $w=\max(0,\sup_{x\in\mathbb{R}}C(x))$. Operator $\overline{H}$ generates a $C_0$-semigroup, so all conditions of the Chernoff theorem are fulfilled. Hence
$$(e^{t\overline{H}}f)(x)=\left(\lim_{n\to\infty}S(t/n)^nf\right)(x), \quad \textrm{for all } \quad t>0, x\in\mathbb{R}, f\in C_0(\mathbb{R},\mathbb{R})$$
locally uniformly in $t$, i.e. for all $T>0$ we have $\lim_{n\to\infty}\sup_{t\in[0,T]}\left\|e^{t\overline{H}}f - (S(t/n))^nf \right\|=0$.

The expression for $S(t/n)^nf$ can be rewritten as follows:
$$
(e^{t\overline{H}}f)(x_0)=\lim_{n\to\infty} \underbrace{\int_{\mathbb{R}}\dots\int_{\mathbb{R}}}_n \exp\left(\frac{t}{n}\sum_{j=1}^n\left(c(x_{j-1})-\frac{b(x_{j-1})^2}{4a(x_{j-1})}\right)\right) \exp\left(\sum_{j=1}^n\frac{b(x_{j-1})(x_j-x_{j-1})}{2a(x_{j-1})} \right)\times
$$
$$
\times\left(\frac{\sqrt{n}}{\sqrt{4\pi t}}\right)^n\left(\prod_{j=0}^{n-1}a(x_j)\right)^{-1/2}\exp\left(-\frac{n}{4t}\sum_{j=0}^{n-1}\frac{(x_j-x_{j+1})^2}{a(x_j)}\right)f(x_n)dx_1\dots dx_n.
$$

\end{remark}

Now we can apply theorem \ref{mainth} and obtain a formula for the resolvent of $\overline{H}$.

\begin{theorem}\label{exth}
Under notation and assumptions from assumption \ref{Hbar} and remark \ref{remhbar}, the resolvent $R_\lambda=(\lambda I-\overline{H})^{-1}$  is given for all $\lambda\in\mathbb{C}$ satisfying $Re\lambda>w$, all $g\in C_0(\mathbb{R},\mathbb{R})$, all $x_0\in \mathbb{R}$ by the following formula:
$$
(R_\lambda g)(x_0) = \lim_{n\to\infty}\int_0^\infty e^{-\lambda t}\Bigg[\underbrace{\int_{\mathbb{R}}\dots\int_{\mathbb{R}}}_n \exp\left(\frac{t}{n}\sum_{j=1}^n\left(c(x_{j-1})-\frac{b(x_{j-1})^2}{4a(x_{j-1})}\right)\right) \exp\left(\sum_{j=1}^n\frac{b(x_{j-1})(x_j-x_{j-1})}{2a(x_{j-1})} \right)\times
$$
$$
\times\left(\frac{\sqrt{n}}{\sqrt{4\pi t}}\right)^n\left(\prod_{j=0}^{n-1}a(x_j)\right)^{-1/2}\exp\left(-\frac{n}{4t}\sum_{j=0}^{n-1}\frac{(x_j-x_{j+1})^2}{a(x_j)}\right)g(x_n)dx_1\dots dx_n\Bigg] dt,
$$
where the limit $\lim\limits_{n\to\infty}$ exists uniformly in $x_0\in\mathbb{R}$.
\end{theorem}
\begin{proof}
Let us check conditions of theorem \ref{mainth}. Set $\mathcal{F}=C_0(\mathbb{R},\mathbb{R})$ with the uniform norm $\|f\|=\sup_{x\in\mathbb{R}}$. Consider $L=\overline{H}$ defined via (\ref{Hlab}), $D(L)=D(\overline{H})$, $\omega=w$, $M=1$. Consider $S(t)$ defined in (\ref{Slab}). We already have $\|S(t)\|\leq e^{wt}$, which implies $\|S(t)^k\|\leq e^{wtk}$. Condition $\lim_{n\to\infty}\sup_{t\in[0,T]}\left\|e^{t\overline{H}}f - (S(t/n))^nf \right\|=0$ for all $f\in\mathcal{F}$ and all $T>0$  is true thanks to \cite{BGS2010}, where the Chernoff theorem is used. The proof is finished due to item 1 of theorem \ref{mainth}.
\end{proof}

\begin{remark}
Note that for an arbitrary linear operator $(L,D(L))$ the resolvent $(\lambda I-L)^{-1}$ and the semigroup $e^{tL}$, if they exist, are defined by $(L,D(L))$ uniquely. Meanwhile,  there are many Chernoff functions for the same operator $(L,D(L))$, so there are many Chernoff approximations for  $(\lambda I-L)^{-1}$ and $e^{tL}$. This gives us some freedom in constructing such approximations. The representation for the resolvent proposed in theorem \ref{exth} is only one of the representations that can be obtained via the Chernoff theorem.
\end{remark}

\section{Corollary: representing solution of ODEs via Feynman formula}

There is no standard well known  method of expressing the solution of ODE  $-a(x)f''(x)-b(x)f'(x)-c(x)f(x)+\lambda f(x)=g(x), x\in\mathbb{R}$ in terms of variable coefficients $a,b,c,g$ and constant $\lambda$. Meanwhile, our method gives a formula for the solution, because $f=(\lambda I-L)^{-1}g$ for $L$ given by $(Lf)(x)=a(x)f''(x)+b(x)f'(x)+c(x)f(x)$. We will rewrite $\lambda f-Lf=g$ as $Lf-\lambda f=-g$ because it is easier to follow the idea. Please allow us to make the statement of the theorem a bit wordy to keep it self-contained. 

\begin{theorem}\label{ODE-th}
Consider second order ordinary differential equation for function $f\colon\mathbb{R}\to\mathbb{R}$ 
\begin{equation}\label{ODEeq}
a(x)f''(x)+b(x)f'(x)+(c(x)-\lambda)f(x)=-g(x)\textrm{ for all }x\in\mathbb{R},
\end{equation}
where functions $a,b,c,g\colon \mathbb{R}\to \mathbb{R}$ are known parameters and number $\lambda\in\mathbb{C}$ is also a known parameter. Assume that there exists constant $a_0>0$ such that $a(x)>a_0$ for all $x\in \mathbb{R}$. Assume that there exists $\beta\in(0,1]$ such that function $c$ is bounded and  H\"older continuous with H\"older exponent $\beta$, and functions $a$, $x\mapsto 1/a(x)$, $b$ are bounded and H\"older continuous with H\"older exponent $\beta$ with derivatives of order one and two. Assume that function $g$ is continuous and vanishes at infinity. Assume that $\mathbb{R}\ni\lambda>\max(0,\sup_{x\in\mathbb{R}}c(x))$.

Then for (\ref{ODEeq}) there exists a unique continuous and vanishing at infinity solution $f$ given for all $x_0\in\mathbb{R}$ by the formula
$$
f(x_0)=\lim_{n\to\infty}\int_0^\infty e^{-\lambda t}\Bigg[\underbrace{\int_{\mathbb{R}}\dots\int_{\mathbb{R}}}_n \exp\left(\frac{t}{n}\sum_{j=1}^n\left(c(x_{j-1})-\frac{b(x_{j-1})^2}{4a(x_{j-1})}\right)\right) \exp\left(\sum_{j=1}^n\frac{b(x_{j-1})(x_j-x_{j-1})}{2a(x_{j-1})} \right)\times
$$
$$
\times\left(\frac{\sqrt{n}}{\sqrt{4\pi t}}\right)^n\left(\prod_{j=0}^{n-1}a(x_j)\right)^{-1/2}\exp\left(-\frac{n}{4t}\sum_{j=0}^{n-1}\frac{(x_j-x_{j+1})^2}{a(x_j)}\right)g(x_n)dx_1\dots dx_n\Bigg] dt,
$$
where the limit $\lim\limits_{n\to\infty}$ exists uniformly in $x_0\in\mathbb{R}$.
\begin{proof} 
In theorem \ref{exth} set $f=R_\lambda g$.
\end{proof}
\end{theorem} 

\begin{remark}
This reasoning also works in the  multi-dimensional situation for  $x\in\mathbb{R}^d$, where we have an elliptic PDE instead of ODE.
\end{remark}

\section{Corollary: translation-based formula as a method of solving ODEs}

Another Chernoff approximations for the same semigroup are known, these approximations do not involve multiple integrals but use multiple shifts instead \cite{RemAMC2018}. For these approximations error bounds are known. Rate of convergence of Chernoff approximations is given in \cite{GR-2021} for the general case of arbitary semigroup, and also in this particular case of translation-based approximations \cite{GR-2021-arxiv} for the semigroup that is discussed in the next theorem. The word ''translation'' is used because for $a(x)=a_0\equiv\mathrm{const}$ operator $f\mapsto [x\mapsto f(x+2\sqrt{a(x)t})]$ is indeed a translation (shift) of $f$ by the value $2\sqrt{a_0t}$.

Let us use symbol $UC_b(\mathbb{R})$ to denote Banach space of all bounded and uniformly continuous functions $f\colon\mathbb{R}\to\mathbb{R}$ with the uniform norm $\|f\|=\sup_{x\in\mathbb{R}}|f(x)|$. Let us use symbol $C_b^\infty(\mathbb{R})$ for the subspace of $UC_b(\mathbb{R})$ consisting of all infinitely differentible functions that are bounded and have bounded derivatives of all orders.

\begin{theorem}\label{thshift}
Suppose that functions $a,b,c\in UC_b(\mathbb{R})$ are bounded with their derivatives up to order 3, and there exists such a constant $a_0>0$ that estimate $\inf_{x\in{\mathbb R}}a(x)\geq a_0>0$ is satisfied for all $x\in\mathbb{R}$.
For each function $\phi\in C_b^\infty(\mathbb{R})=D(A)$ define $A\phi = a\phi'' + b\phi' + c\phi$. For each $t\geq0$, each $x\in\mathbb{R}$ and each $f\in UC_b(\mathbb{R})$ define 
\begin{equation}\label{Sdef}
			(S(t)f)(x) = \frac14 f\Big(x+2\sqrt{a(x)t}\Big) + \frac14 f\Big(x-2\sqrt{a(x)t}\Big) 
			+ \frac12 f\big(x+2b(x)t\big) + tc(x)f(x).
	\end{equation}

Assume also that $\mathbb{R}\ni\lambda>\sup_{x\in\mathbb{R}}|c(x)|=\|c\|$. Then: 

1. Closure $\overline{A}$ of operator $A$ generates a $C_0$-semigroup in $UC_b(\mathbb{R})$.

2. For each $g\in UC_b(\mathbb{R})$ the solution $f\colon\mathbb{R}\to\mathbb{R}$ of the equation
\begin{equation}\label{ODEeq2}
	a(x)f''(x)+b(x)f'(x)+(c(x)-\lambda)f(x)=-g(x)\textrm{ for all }x\in\mathbb{R},
\end{equation}
exists, is unique in $UC_b(\mathbb{R})$ and is given for all $x\in\mathbb{R}$ by the formula
\begin{equation}\label{ODEsol}
f(x)=\int_0^\infty e^{-\lambda t}\left(e^{\overline{A}}g\right)(x)dt=\lim_{n\to\infty}\int_0^\infty e^{-\lambda t}\left((S(t/n))^ng\right)(x)dt,
\end{equation}
where $S(t/n)$ is obtained by substitution of $t$ with $t/n$ in (\ref{Sdef}), and $(S(t/n))^n$ is the composition of $n$ copies of linear bounded operator $S(t/n)$.

Suppose additionally that function $g$ is bounded with derivatives up to order 5. Then:

3. There exist nonnegative constants $C_0,C_1,\ldots,C_4$ such that 
	for all $t>0$ and all $n\in\mathbb{N}$ the following inequality holds:
	\begin{equation} \| S(t/n)^n g - e^{t\overline{A}}g \| \leq \frac{t^2e^{\|c\|t}}{n}\big(C_0\|g\|+C_1\|g'\|+C_2\|g''\|+C_3\|g'''\|+C_4\|g^{(IV)}\|\big).
	\end{equation}

4. Error bound in (\ref{ODEsol}) for all $n\in\mathbb{N}$ is given by inequality
$$
\sup_{x\in\mathbb{R}}\left|f(x)-\int_0^\infty e^{-\lambda t}\left((S(t/n))^ng\right)(x)dt\right|\leq\frac{2C_g}{n\cdot (\lambda-\|c\|)^3},
$$
where $C_g=C_0\|g\|+C_1\|g'\|+C_2\|g''\|+C_3\|g'''\|+C_4\|g^{(IV)}\|$.

5. Integral in item 2 can be calculated over $[0,T]$ instead of $[0,\infty)$ with controlled level of error. This means that for each $\varepsilon>0$ there exists $
T=\max\left(0,\frac{1}{\lambda-\|c\|}\ln\frac{2}{(\lambda-\|c\|)\varepsilon}\right)$ such that for all $n\in\mathbb{N}$ we have
$$
\sup_{x\in\mathbb{R}}\left|f(x)-\int_0^Te^{-\lambda t}\left((S(t/n))^ng\right)(x)dt\right|\leq\frac{2C_g}{n\cdot (\lambda-\|c\|)^3}+\varepsilon.
$$
\end{theorem}

\begin{proof}
Item 1 follows from theorem 4.2 in \cite{GR-2021-arxiv}. Item 2 is a particular case of the main result of the paper, theorem \ref{mainth}. Item 3 follows is example 4.2 in \cite{GR-2021-arxiv}. Item 4 follows from items 2 and 3 with simple estimate
\begin{multline*}
\left\|\int_0^\infty e^{-\lambda t}\left(e^{\overline{A}}g\right)(x)dt-\int_0^\infty e^{-\lambda t}\left((S(t/n))^ng\right)(x)dt\right\|\leq\\ 
\leq\int_0^\infty e^{-\lambda t}\left\|e^{\overline{A}}g-(S(t/n))^ng\right\|dt\leq \int_0^\infty e^{(\|c\|-\lambda)t}\frac{t^2}{n}C_gdt=\frac{2C_g}{n\cdot (\lambda-\|c\|)^3}.
\end{multline*}	
Item 5 (by repeating the reasoning in the first part of the proof of theorem \ref{mainth}) folows from item 4 and the well known fact that the semigroup $(e^{t\overline{A}})_{t\geq0}$ is a quasi-contraction, i.e. in estimate for norm $\|e^{t\overline{A}}\|\leq Me^{\omega t}$ it is possible to set $M=1$, $\omega=\|c\|$.
\end{proof}

\begin{remark}
Independently of Chernoff function used (is it based on integral operators as in theorem \ref{ODE-th} or on translation operators as in theorem \ref{thshift}), Chernoff approximations are allowing to calculate value of the solution in only one point of the domain of solution (in one point $x\in\mathbb{R}$ in our examples). Meanwhile methods based on a computational grid calculate values of the solution in all points of the computational grid. Moreover, values of Chernoff approximations at different points of the domain can be calculated in parallel, using multi-core processors and GPU which is an advantage of this approach. 
\end{remark}

\textbf{Acknowledgments.} Author is thankful to Oleg Galkin, Denis Mineev and Polina Panteleeva for comments on the manuscript. The work (except theorem 5) was supported by the Russian Science Foundation (project 23-71-30008). Theorem 5 was obtained in IITP (Dobrushin's Math. Lab.).


\end{document}